\documentclass[12pt]{amsart}
\usepackage[all]{xy}
\usepackage{a4wide}
\usepackage{mathtools}
\usepackage{amssymb}
\usepackage{amsthm}
\usepackage{amsmath}
\usepackage{amscd}
\usepackage{verbatim}
\usepackage[all]{xy}
\usepackage{xcolor}

\usepackage{amsmath,amsthm,amssymb}
\usepackage{enumitem}

\usepackage{float}
\usepackage{ifthen}
\usepackage{xargs}
\usepackage{url}
\usepackage[all]{xy}
\usepackage{xcolor}
\usepackage{comment}
\usepackage[obeyDraft,textsize=tiny]{todonotes}
\specialcomment{draftnote}{\vspace{1em}\begingroup\ttfamily\footnotesize\color{blue}\begin{minipage}{15cm}}{\end{minipage}\endgroup\vspace{1em}}
\specialcomment{alert}{\vspace{1em}\begingroup\ttfamily\footnotesize\color{blue}\begin{minipage}{15cm}}{\end{minipage}\endgroup\vspace{1em}}

\usepackage[refpage,intoc,noprefix]{nomencl}

\newcommand{\R}{\mathbb{R}} 
\newcommand{\Q}{\mathbb{Q}} 
\newcommand{\C}{\mathbb{C}} 
\newcommand{\Z}{\mathbb{Z}} 
\newcommand{\N}{\mathbb{N}} 
\renewcommand{\H}{\mathbb{H}} 

\newcommand*{\reg}{\mathrm{reg}}

\newcommand{\tr}{\mathrm{tr}}

\newcommand{\SL}{{\text {\rm SL}}}

\newcommand{\e}{\mathfrak{e}}

\DeclareMathOperator{\Mp}{Mp}

\DeclareMathOperator{\Gr}{Gr}

\DeclareMathOperator{\sgn}{sgn}

\newtheorem{theorem}{Theorem}[section]
\newtheorem{proposition}[theorem]{Proposition}
\newtheorem{lemma}[theorem]{Lemma}

\theoremstyle{definition}

\newtheorem{remark}[theorem]{Remark}

\numberwithin{theorem}{section} \numberwithin{equation}{section}

\newcommand{\new}{\mathrm{new}}

\definecolor{ocre}{rgb}{0.37, 0.62, 0.63} 
\definecolor{redocre}{RGB}{156, 95, 94} 
\definecolor{greenocre}{RGB}{95, 155, 95} 
\definecolor{darkblueocre}{RGB}{95, 94, 156} 
\definecolor{orangeocre}{RGB}{184, 142, 71} 
\usepackage[framemethod=tikz]{mdframed}
\mdtheorem[
  linecolor=darkblueocre,
  frametitlefont=\sffamily\bfseries\color{white},
  frametitlebackgroundcolor=darkblueocre,
]{boxthm}{Theorem}[section]
\mdtheorem[
  linecolor=darkblueocre,
  frametitlefont=\sffamily\bfseries\color{white},
  frametitlebackgroundcolor=darkblueocre,
]{boxprop}{Proposition}[section]

\mdtheorem[
  linecolor=ocre,
  frametitlefont=\sffamily\bfseries\color{white},
  frametitlebackgroundcolor=ocre,
]{boxthmold}{Theorem}[section]

\begin{document}

\title[]{Harmonic weak Maass forms and periods II}

\author{Claudia Alfes-Neumann, Jan Hendrik Bruinier, Markus Schwagenscheidt}

\begin{abstract}
	In this paper we investigate the Fourier coefficients of harmonic Maass forms of negative half-integral weight. We relate the algebraicity of these coefficients to the algebraicity of the coefficients of certain canonical meromorphic modular forms of positive even weight with poles at Heegner divisors. Moreover, we give an explicit formula for the coefficients of harmonic Maass forms in terms of periods of certain meromorphic modular forms with algebraic coefficients.
\end{abstract}

\maketitle

\section{Introduction}

A fundamental result in the theory of modular forms is the fact that the spaces of holomorphic modular forms of fixed weight possess bases of forms with \emph{integral} Fourier coefficients. In recent years, many authors have studied the Fourier coefficients of non-holomorphic generalizations of modular forms, such as \emph{harmonic Maass forms}. These functions transform like modular forms but are harmonic rather than holomorphic on the upper half-plane, and they are allowed to have poles at the cusps. The most prominent examples appear in connection with mock theta functions, which are certain $q$-series with (typically) integral coefficients that have asymptotic expansions like modular forms as $q$ approaches roots of unity, but are not modular. More precisely, Zwegers \cite{zwegers} showed that Ramanujan's mock theta functions may be viewed as the \emph{holomorphic parts} of certain harmonic Maass forms of weight $\frac{1}{2}$. This raises the question for the algebraicity of the coefficients of the holomorphic parts of more general harmonic Maass forms. 

As it turns out, apart from the mock theta functions (see also \cite{brfolono,bs}) and certain harmonic Maass forms related to newforms with complex multiplication (see \cite{bor}), these coefficients typically do not seem to be algebraic. For example, in \cite{bruinieronoheegner}, Ono and the second author related the algebraicity of the Fourier coefficients of the holomorphic parts of certain harmonic Maass forms of weight $\frac{1}{2}$ to the vanishing of the central values of derivatives of $L$-functions of newforms of weight $2$. Moreover, in \cite{brdif} it was shown that the Fourier coefficients of these weight $\frac{1}{2}$ harmonic Maass forms can be expressed in terms of periods of certain algebraic differentials of the third kind on modular curves, which also hints at their delicate algebraic nature. In the present work, we extend the results of \cite{brdif} to harmonic Maass forms of weight $\frac{3}{2}-k$ with $k \geq 1$, and derive an expression for their Fourier coefficients in terms of periods of certain meromorphic modular forms of weight $2k$.

Let us describe our main results in some more detail. To simplify the exposition we let $\Gamma:=\SL_2(\Z)$ throughout the introduction. In the body of the paper we will treat level $\Gamma_0(N)$ for arbitrary $N \in \N$, using the language of vector-valued modular forms for the Weil representation.

\subsection{Canonical meromorphic modular forms associated with divisors on modular curves} If $X$ is a non-singular projective curve over $\C$ and $D$ is a divisor on $X$ whose restriction to any component of $X$ has degree $0$, then it follows from the Riemann period relations that there exists a unique meromorphic differential $\eta_D$ on $X$ with at most simple poles and residue divisor $\sum_{P \in X}\mathrm{res}_P(\eta_D)P = D$, and such that
\[
\Re\left(\int_C \eta_D \right) = 0
\]
for any closed path $C$ in $X$ avoiding the points of $D$. The differential $\eta_D$ is called the \emph{canonical differential of the third kind associated with $D$}. If $X$ is a modular curve, then $\eta_D$ may be viewed as a meromorphic modular form of weight $2$. Our first result is a construction of certain \lq canonical\rq \ meromorphic modular forms $\eta_{k,D}$ of weight $2k$ associated with divisors $D$ on modular curves, generalizing the notion of canonical differentials to higher weight. 

For a subfield $F \subseteq \C$ and an integer $k \in \Z$ with $k \geq 2$ we let $\mathcal{D}_{2k,F}(\Gamma)$ be the $F$-vector space of meromorphic modular forms $f(z)$ of weight $2k$ for $\Gamma$ which vanish at the cusp of $\Gamma$ and which have poles of order $k$ on $\H$, with expansions of the form
\[
f(z) = a_{f,\varrho}\left(\frac{(z-\varrho)(z-\overline{\varrho})}{\varrho-\overline{\varrho}}\right)^{-k}+ \mathcal{O}(1), \qquad (z \to \varrho),
\]
at each pole $\varrho$ of $f$, with some constants $a_{f,\varrho} \in F$. Note that for $k = 1$ we have $a_{f,\varrho} = \mathrm{res}_{\varrho}(f)$. Hence, by a slight abuse of notation, we will call the $F$-linear combination
\[
\mathrm{res}(f) := \sum_{[\varrho] \in \Gamma\backslash \H}\frac{a_{f,\varrho}}{w_\varrho}[\varrho] \in \mathrm{Div}(\Gamma\backslash \H)_F
\]
of points in $\Gamma\backslash\H$ the \emph{residue divisor} of $f \in \mathcal{D}_{2k,F}(\Gamma)$.  Here $w_\varrho:= \frac{1}{2}|\Gamma_\varrho|$ is half the order of the stabilizer of $\varrho$ in $\Gamma$, and $\mathrm{Div(\Gamma\backslash \H)}_{F}$ denotes the space of all $F$-linear combinations of points on $\Gamma \backslash \H$. Note that the sum is well-defined, that is, $a_{f,\varrho}$ and $w_\varrho$ only depend on the class $[\varrho]$ of $\varrho$ in $\Gamma\backslash \H$.

It is well known that closed geodesics on $\Gamma \backslash \H$ correspond to (conjugacy classes of) primitive hyperbolic matrices in $\Gamma$. More precisely, given such a primitive hyperbolic matrix $\gamma \in \Gamma$, we let $S_\gamma$ be the semi-circle in $\H$ connecting the two real fixed points of $\Gamma$. Then $C_\gamma := \langle \gamma \rangle \backslash S_\gamma$ defines a closed geodesic in $\Gamma \backslash \H$, and every closed geodesic in $\Gamma \backslash \H$ is accounted for in this way. We let $\mathcal{C}_{\R}(\Gamma)$ be the real vector space spanned by all $\R$-linear combinations of closed geodesics $C_\gamma$ on $\Gamma \backslash \H$. We obtain a real-valued bilinear pairing on $\mathcal{D}_{2k,\R}(\Gamma) \times \mathcal{C}_{\R}(\Gamma)$ by setting, for each meromorphic modular form $f \in \mathcal{D}_{2k,\R}(\Gamma)$ and each closed geodesics $C_\gamma$ on $\Gamma \backslash \H$,
\begin{align}\label{pairing}
(f,C_\gamma) := \Re\left(\left(\frac{i}{\sqrt{d_\gamma}}\right)^{k-1}\int_{z_0}^{\gamma z_0}f(z)Q_\gamma(z,1)^{k-1}dz\right),
\end{align}
where $Q_\gamma(x,y) := cx^2 + (d-a)xy - by^2$ denotes the binary quadratic form corresponding to $\gamma = \left(\begin{smallmatrix}a & b \\ c & d \end{smallmatrix}\right)$, $d_\gamma := \tr(\gamma)^2-4$ is the discriminant of $Q_\gamma$, and the path of integration is any path in $\H$ (avoiding the poles of $f$) from some point $z_0 \in \H$ (not being a pole of $f$) to $\gamma z_0$. Using the residue theorem and the assumption that the coefficients $a_{f,\varrho}$ are real numbers, one can show that the value of the pairing $(f,C_\gamma)$ is indeed independent of the choice of $z_0 \in \H$ and the path from $z_0$ to $\gamma z_0$, see Lemma~\ref{pairing path independent} below. 

We have the following higher weight analogues of canonical differentials associated with divisors on $\Gamma \backslash \H$. 

\begin{proposition}\label{canonical modular form}
	For each divisor $D \in \mathrm{Div}(\Gamma\backslash \H)_{\R}$ there exists a unique meromorphic modular form $\eta_{k,D} \in \mathcal{D}_{2k,\R}(\Gamma)$ with residue divisor $\mathrm{res}(\eta_{k,D}) = D$ and $(\eta_{k,D},C) = 0$ for any $C \in \mathcal{C}_{\R}(\Gamma)$.
\end{proposition}

We call $\eta_{k,D}$ the \emph{canonical meromorphic modular form of weight $2k$ associated with $D$}. We will construct it explicitly as a linear combination of certain two-variable Poincar\'e series introduced by Petersson \cite{petersson}. The fact that their cycle integrals are all real or all purely imaginary (depending on the parity of $k$) can then be shown by a direct computation, following ideas of Katok \cite{katok}. Moreover, the uniqueness of $\eta_{k,D}$ essentially follows from the Eichler-Shimura isomorphism. We refer the reader to Section~\ref{meromorphic modular forms} for the details of the proof of Proposition~\ref{canonical modular form}.

In the case that $D$ is a CM divisor, a result similar to Proposition~\ref{canonical modular form} has been proved by Mellit \cite[Theorem~1.5.3]{mellit}.

\begin{remark}
The above constructions may be interpreted in a more geometric way using local coefficient systems as explained in \cite{funkemillson}. In this setting, meromorphic modular forms of weight $2k$ give rise to meromorphic $1$-forms on $\Gamma \backslash \H$ with coefficients in the local system $E$ associated with $\mathrm{Sym}^{2k-2}(\C^2)$, geodesics in $\Gamma\backslash \H$ give rise to cycles with coefficients in the dual of $E$, and the pairing \eqref{pairing} may be viewed as a (co)homological pairing.
\end{remark}

\subsection{Harmonic Maass forms and normalized meromorphic modular forms} 

The main topic of this work is the study of canonical meromorphic modular forms of weight $2k$ for certain linear combinations of \emph{Heegner divisors}, coming from harmonic Maass forms of half-integral weight. Let us recall the setup from \cite{bruinieronoheegner}. For $k \in \N$ with $k \geq 2$ we let $G \in S_{2k}$ be a normalized Hecke eigenform of weight $2k$ for $\Gamma$, and let $F_G$ be the totally real number field generated by the Fourier coefficients of $G$. We let $g \in S_{\frac{1}{2}+k}$ be a Hecke eigenform of weight $\frac{1}{2}+k$ for $\Gamma_0(4)$ in the Kohnen plus space which corresponds to $G$ under the Shimura correspondence. We may normalize $g$ such that all its Fourier coefficients lie in $F_G$. Then, by \cite[Lemma~7.3]{bruinieronoheegner} there exists a harmonic Maass form $f$ of weight $\frac{3}{2}-k$ for $\Gamma_0(4)$ with principal part defined over $F_G$ such that $\xi_{\frac{3}{2}-k}(f) = \|g\|^{-2}g$, where $\|g\| :=\sqrt{(g,g)}$ denotes the Petersson norm of $g$ and $\xi_{\kappa} = \overline{2i v^{\kappa}\frac{\partial}{\partial \overline{\tau}}}$. Every such harmonic Maass form $f$ can be written as a sum $f=f^+ + f^-$ of a \emph{holomorphic part} $f^+$ and a \emph{non-holomorphic part} $f^-$, with Fourier expansions of the shape
\[
f^+(\tau) =\sum_{\substack{n \in \Z \\ n \gg -\infty}}c_f^+(n)q^n, \qquad f^-(\tau) = \sum_{\substack{n \in \Z \\ n < 0}}c_f^-(n)\Gamma\left(k-\tfrac{1}{2},4\pi|n|v\right)q^n,
\] 
with coefficients $c_f^{\pm}(n) \in \C$, where $q := e^{2\pi i\tau}$ and $\Gamma(s,x) := \int_{x}^\infty e^{-t}t^{s-1}dt$ denotes the incomplete Gamma function. The condition that the principal part of $f$ is defined over $F_G$ means that we have $c_f^{+}(n) \in F_G$ for all $n \leq 0$. We let $H_{\frac{3}{2}-k}(F_G)$ denote the space of harmonic Maass forms of weight $\frac{3}{2}-k$ for $\Gamma_0(4)$ satisfying the Kohnen plus space condition and with principal part defined over $F_G$. We are interested in the algebraicity properties of the coefficients $c^+_f(n)$ for $n > 0$.

Now we come to the construction of the Heegner divisor corresponding to $f$. For a negative discriminant $d < 0$ we let $\mathcal{Q}_d^+$ be the set of positive definite integral binary quadratic forms $Q(x,y) = ax^2+bxy +cy^2$ of discriminant $d = b^2-4ac$. For $Q \in \mathcal{Q}_d^+$ the quadratic polynomial $Q(z,1)$ has a unique root $\alpha_Q$ in $\H$, called the \emph{CM point} (or \emph{Heegner point}) corresponding to $Q$. Let $\Delta \in \Z$ be a fundamental discriminant with $(-1)^k \Delta < 0$ and let $D < 0$ be an integer such that $|\Delta|D \equiv 0,1 \pmod 4$. Moreover, let $\chi_\Delta$ be the usual genus character on $\mathcal{Q}_{D|\Delta|}^+$. Then we define the \emph{twisted Heegner divisor}
\[
Z_{\Delta}(D) := \sum_{Q \in \Gamma \backslash\mathcal{Q}_{D|\Delta|}^+}\frac{\chi_{\Delta}(Q)}{w_Q}[\alpha_Q],
\]
where $w_Q := w_{\alpha_Q}$. Eventually, we let
\[
Z_{k,\Delta}(f) := \sum_{D < 0}c_f^+(D)\big(|D\Delta|\big)^{\frac{k-1}{2}}Z_{\Delta}(D)
\]
be the twisted Heegner divisor corresponding to $f$, where $c_f^+(D) \in F_G$ denote the coefficients of the holomorphic part of $f$. We remark that for even $k$ the Heegner divisor $Z_{k,\Delta}(f)$ is in general not a divisor with coefficients in $F_G$, but with coefficients in some totally real algebraic extension of $F_G$.

 The following result relates the algebraicity of the coefficients $c_f^+(n)$ of the harmonic Maass form $f$ to the algebraicity of the Fourier coefficients of the canonical meromorphic modular form of weight $2k$ for the Heegner divisor $Z_{k,\Delta}(f)$. It (partly) generalizes \cite[Theorem~5.5, Theorem~7.6]{bruinieronoheegner} to higher weight.

\begin{theorem}\label{canonical modular form heegner}
	Let $f \in H_{\frac{3}{2}-k}(F_G)$ be as above, and let $\eta_{k,\Delta}(f) \in \mathcal{D}_{2k,\R}(\Gamma)$ be the canonical meromorphic modular form of weight $2k$ for $Z_{k,\Delta}(f)$. Then the following are equivalent.
	\begin{enumerate}
		\item We have $c_{f}^+(|\Delta|) \in F_G$.
		\item We have $c_{f}^+(n^2|\Delta|) \in F_G$ for all $n \in \N$.
		\item All Fourier coefficients of $\eta_{k,\Delta}(f)$ are contained in $i\pi^k \sqrt{\Delta} F_G$.
	\end{enumerate}
 \end{theorem}

 In order to prove the theorem we will construct $\eta_{k,\Delta}(f)$ as a regularized theta lift of $f$, extending the methods of \cite{bruinieronoheegner} by including certain differential (weight raising) operators in the theta lift, similarly as in \cite{bruinierehlenyang}. By computing the Fourier expansion of this lift, we obtain a formula for the coefficients of $\eta_{k,\Delta}(f)$ in terms of the coefficients $c_f^+(n^2|\Delta|), n \in \N$. Invoking the action of Hecke operators on $f$ and $\eta_{k,\Delta}(f)$ then yields the stated result. We refer to Section~\ref{meromorphic modular forms for heegner divisors} for the generalization of Theorem~\ref{canonical modular form heegner} to higher level $\Gamma_0(N)$ and its proof.
 
 We remark that a similar theta lift has been studied by Zemel~\cite{zemel}, and it was used to show the implication $(2)\Rightarrow(3)$ in Theorem~\ref{canonical modular form heegner} (for $\Delta = 1$).

Finally, we obtain a formula for the coefficient $c_f^+(|\Delta|)$ in terms of cycle integrals of a certain \emph{normalized} meromorphic modular form in $\mathcal{D}_{2k,\R}(\Gamma)$ with residue divisor $Z_{k,\Delta}(f)$, generalizing \cite[Theorem~1.1]{brdif} to higher weight.

\begin{theorem}\label{normalized modular form heegner}
	For $f \in H_{\frac{3}{2}-k}(F_G)$ as above, there exists a unique meromorphic modular form $\zeta_{k,\Delta}(f) \in \mathcal{D}_{2k,\R}(\Gamma)$ with the following properties.
	\begin{enumerate}
		\item The residue divisor is given by $\mathrm{res}(\zeta_{k,\Delta}(f)) = Z_{k,\Delta}(f)$.
		\item We have $(\zeta_{k,\Delta}(f),C) = 0$ for every $C \in \mathcal{C}_{\R}(\Gamma)$ with $(G,C) = 0$.			
		\item The first Fourier coefficient of $\zeta_{k,\Delta}(f)$ vanishes.
	\end{enumerate}
	Moreover, all Fourier coefficients of $\zeta_{k,\Delta}(f)$ lie in $i\pi^k \sqrt{\Delta}F_G$, and we have
	\begin{align}\label{period formula intro}
	c_f^+(|\Delta|) = -\frac{1}{C_{k,\Delta}i\pi^k\sqrt{\Delta}}\cdot\frac{(\zeta_{k,\Delta}(f),C)}{(G,C)}
	\end{align}
	for every $C \in \mathcal{C}_{\R}(\Gamma)$ with $(G,C) \neq 0$, with the rational constant $C_{k,\Delta} = \frac{(-2)^{k}|\Delta|^{k-1}}{(k-1)!}$.
\end{theorem}

 \begin{remark}
 The formula from Theorem~\ref{normalized modular form heegner} tells us that the coefficient $c_f^+(|\Delta|)$ can be written as a quotient of (real parts of) cycle integrals of (meromorphic) modular forms with coefficients in $F_G$.
 \end{remark}

We call $\zeta_{k,\Delta}(f)$ the \emph{normalized meromorphic modular form of weight $2k$ for $Z_{k,\Delta}(f)$}.  We will construct $\zeta_{k,\Delta}(f)$ by subtracting a suitable multiple of $G$ from $\eta_{k,\Delta}(f)$. The fact that the coefficients $\zeta_{k,\Delta}(f)$ lie in $i\pi^k \sqrt{\Delta}F_G$ and the formula \eqref{period formula intro} again follow by writing $\eta_{k,\Delta}(f)$ as a theta lift of $f$, and using the action of Hecke operators on $\zeta_{k,\Delta}(f)$. The details of the proof can be found in Section~\ref{meromorphic modular forms for heegner divisors}.

 \subsection{Organization of the paper} We start with a section on the necessary preliminaries about Heegner divisors, quadratic spaces and lattices, and vector-valued harmonic Maass forms for the Weil representation. In Section~\ref{theta lifts} we study a regularized theta lift of harmonic Maass forms which produces canonical meromorphic modular forms for twisted Heegner divisors, and we compute its Fourier expansion. In Section~\ref{meromorphic modular forms for heegner divisors} we investigate the algebraicity properties of the Fourier coefficients of the canonical and normalized meromorphic modular forms for twisted Heegner divisors, and we prove Theorem~\ref{canonical modular form heegner} and Theorem~\ref{normalized modular form heegner} for level $\Gamma_0(N)$. Finally, in Section~\ref{outlook} we give a short outlook on a possible future application of our results to a non-vanishing criterion for central values of derivatives of newforms of weight $2k$ in terms of the algebraicity of Fourier coefficient of harmonic Maass forms.
 
 \subsection*{Acknowledgment} We thank Jens Funke for helpful discussions. The first author was partially supported by the Daimler and Benz Foundation and the Klaus Tschira Boost Fund. The second author was supported by the DFG Collaborative Research Centre TRR~326 {\em Geometry and Arithmetic of Uniformized Structures}, project number 444845124. The third author was supported by SNF projects 200021\_185014 and PZ00P2\_202210.

\section{Preliminaries}
\label{preliminaries}

\subsection{Heegner divisors} 
\label{quadratic forms}

Let $N\in \N$. For a negative discriminant $d < 0$ and $r \in \Z/2N\Z$ with $d \equiv r^2 \pmod{4N}$ we consider the set of integral binary quadratic forms
\[
\mathcal{Q}_{d,r} := \{Q(x,y) = aNx^2 + bxy+cy^2 : a,b,c \in \Z \, , \, d = b^2-4Nac \, ,\, b\equiv r \!\!\!\!\pmod{2N}\}.
\]
The group $\Gamma_0(N)$ acts on $\mathcal{Q}_{d,r}$ from the right, with finitely many orbits if $D \neq 0$.

We let $\mathcal{Q}_{d,r}^+$ be the subset of positive definite quadratic forms in $\mathcal{Q}_{d,r}$, which are given by the condition $\sgn(Q):=\sgn(a) > 0$. For $Q = [a,b,c] \in \mathcal{Q}_{d,r}^+$ we let 
\[
\alpha_Q := \frac{-b+i\sqrt{|d|}}{2Na}
\]
be the \emph{Heegner point} corresponding to $Q$. It is the unique root of $Q(z,1)$ in $\H$. The order $w_Q := \frac{1}{2}|\Gamma_0(N)_Q|$ of the stabilizer of $Q$ in $\Gamma_0(N)$ is finite.

Let $\Delta \in \Z$ be a fundamental discriminant and let $\rho \in \Z/2N\Z$ with $\Delta \equiv \rho^2 \pmod{4N}$. For $D \in \Z$ with $D < 0$ we let $\chi_\Delta$ be the generalized genus character on $\mathcal{Q}_{D|\Delta|,r\rho}$ as in \cite{gkz}, and we define the \emph{twisted Heegner divisor}
\[
Z_{\Delta,\rho}(D,r) := \sum_{Q \in \mathcal{Q}_{D|\Delta|,r\rho}^+/\Gamma_0(N)}\frac{\chi_\Delta(Q)}{w_Q}[\alpha_Q],
\]
where $[\alpha_Q]$ denotes the class of $\alpha_Q$ in $\Gamma_0(N)\backslash \H$.

\subsection{Quadratic spaces and lattices}
\label{quadratic spaces}

For $N \in \N$ we let $V$ be the rational quadratic space of signature $(2,1)$ given by the set of rational traceless $2$ by $2$ matrices with the quadratic form $q(X) := -N\det(X)$ and the associated bilinear form $(X,Y) := N\tr(XY)$. The group $\SL_2(\Q)$ acts as isometries on $V$ by $\gamma.X := \gamma X\gamma^{-1}$.

 Throughout this work we will consider the lattice
\begin{align}\label{the lattice}
L :=  \left\{\begin{pmatrix}b & c/N \\ -a & -b\end{pmatrix}: a,b,c \in \Z \right\}
\end{align}
in $V$. The dual lattice is given by
\[
L' = \left\{\begin{pmatrix}b/2N & c/N \\ -a & -b/2N\end{pmatrix}: a,b,c \in \Z \right\}.
\]
Hence the discriminant form $L'/L$ is isomorphic to $\Z/2N\Z$ with the finite quadratic form $x \mapsto x^2/4N \pmod \Z$, and we will use this identification without further notice. Note that the group $\Gamma_0(N)$ acts on $L$ and fixes the classes of $L'/L$.

For $r \in \Z/2N\Z$ with $d \equiv r^2 \pmod{4N}$ we let 
\[
L_{d,r} := \{X \in L': q(X) = d/4N \text{ and }X \equiv r \!\!\!\! \pmod L\},
\]
Each element $X = \left(\begin{smallmatrix}b/2N & c/N \\ -a & -b/2N \end{smallmatrix}\right)  \in L_{d,r}$ corresponds to a binary quadratic form 
\[
Q_X := [aN,b,c] \in \mathcal{Q}_{d,r}
\]
of discriminant $d = 4Nq(X)$. The group $\Gamma_0(N)$ acts on both $L_{d,r}$ and $\mathcal{Q}_{d,r}$, and the actions are compatible in the sense that $Q_X \circ \gamma = Q_{\gamma^{-1}.X}$ for $\gamma \in \Gamma_0(N)$. In particular, we obtain a bijection between $\Gamma_0(N)\backslash L_{d,r}$ and $\mathcal{Q}_{d,r}/\Gamma_0(N)$.

We let $\Gr(L)$ be the Grassmannian of positive definite $2$-dimensional subspaces of $V(\R)$. We can identify $\Gr(L)$ with $\H$ by sending $z = x+iy \in \H$ to the positive definite plane $U(z) := X(z)^\perp$, where
\[
X(z) := \frac{1}{\sqrt{2N}y}\begin{pmatrix} -x & |z|^2 \\ -1 &x \end{pmatrix}.
\]
Note that $(X(z),X(z)) = -1$ and $\gamma.X(z) = X(\gamma z)$ for $\gamma \in \SL_2(\R)$. We also define the vectors
\begin{align*}
U_1(z) := \frac{1}{\sqrt{2N}y}\begin{pmatrix}x & -x^2+y^2 \\ 1 & -x \end{pmatrix},\qquad
U_2(z) := \frac{1}{\sqrt{2N}y}\begin{pmatrix}y & -2xy \\ 0 & -y \end{pmatrix},
\end{align*}
such that $U_1(z),U_2(z),X(z)$ form an orthogonal basis of $L \otimes \R$. Note that $U(z)$ is spanned by $U_1(z),U_2(z)$. For $X \in L'$ and $z \in \H$ we consider the polynomials
\begin{align}
\begin{split}\label{polynomials}
Q_X(z) &:= \sqrt{2N}y(X,U_1(X) + iU_2(X)) = aNz^2 + bz + c, \\
p_z(X) &:= -(X,X(z)) = \frac{1}{\sqrt{2N}y}(aN|z|^2 + bx + c).
\end{split}
\end{align}
Then we have the useful rules
\begin{align}
\label{pq identity}
q(X_z) = \frac{1}{4Ny^2}|Q_X(z)|^2,\qquad 
q(X_{z^\perp}) &= -\frac{1}{2}p_z^2(X).
\end{align}

\subsection{Harmonic Maass forms}
\label{harmonic maass forms}

The {\it metaplectic extension of $\SL_2(\Z)$} is defined as
\begin{equation*}
\widetilde{\Gamma} \coloneqq \text{Mp}_2(\Z) \coloneqq \left\{ (\gamma, \phi) \colon \gamma = \left(\begin{matrix}
a & b \\ c & d
\end{matrix}\right)\in \SL_2(\Z), \phi\colon \H \rightarrow \C \text{ holomorphic}, \phi^2(\tau) = c\tau+d  \right\}.
\end{equation*}
It is generated by $T := \left(\left( \begin{smallmatrix} 1 & 1 \\ 0 & 1 \end{smallmatrix} \right),1\right)$ and $S := \left(\left( \begin{smallmatrix} 0 & -1 \\ 1 & 0 \end{smallmatrix}\right) ,\sqrt{\tau}\right)$. We let $\widetilde{\Gamma}_\infty$ be the subgroup generated by $T$. Moreover, for $\gamma = \left(\begin{smallmatrix} a & b \\ c & d \end{smallmatrix} \right) \in \SL_2(\Z)$ we let $\widetilde{\gamma} := (\gamma,\sqrt{c\tau + d}) \in \widetilde{\Gamma}$, where $\sqrt{\cdot}$ denotes the principal branch of the square root.

We let $L$ be the lattice from Section~\ref{quadratic spaces} and we let $\C[L'/L] \cong \C[\Z/2N\Z]$ be its group ring, which is generated by the formal basis vectors $\e_r$ for $r \in\Z/2N\Z$. The {\it Weil representation} $\rho_L$ associated with $L$ is the representation of $\widetilde{\Gamma}$ on $\C[L'/L]$ defined by
\begin{align*}
\rho_L(T)(\e_r) \coloneqq e\left(\frac{r^2}{4N}\right) \e_r, \qquad
\rho_L(S)(\e_r) \coloneqq \frac{e\left(-\frac18\right)}{\sqrt{2N}} \sum_{r'(2N)} e\left(-\frac{rr'}{2N}\right) \e_{r'}.
\end{align*}
The Weil representation $\overline{\rho}_{L}$ associated to the lattice $L^- = (L,-q)$ is called the \textit{dual Weil representation} associated to $L$.

Let $\kappa \in \Z+\frac{1}{2}$ and define the {\it slash-operator} by
\[
f\mid_{\kappa,\rho_{L}}(\gamma,\phi) (\tau):= \phi(\tau)^{-2\kappa}\rho_{L}^{-1}(\gamma,\phi)f(\gamma\tau),
\]
for a function $f\colon \H \rightarrow \C[L'/L]$ and $(\gamma,\phi) \in \widetilde{\Gamma}$. Here, and in the following we let $e(z):=e^{2\pi i z}$. Following \cite{brfu04}, we call a smooth function $f \colon \H \rightarrow \C[L'/L]$ a {\it harmonic  Maass form} of weight $\kappa$ with respect to $\rho_L$ if it is annihilated by the \textit{weight $\kappa$} \textit{Laplace operator}
\[
\Delta_{\kappa} \coloneqq -v^2\left(\frac{\partial^2}{\partial u^2} + \frac{\partial^2}{\partial v^2} \right) + i\kappa v\left(\frac{\partial}{\partial u} + i\frac{\partial}{\partial v} \right),
\]
if it is invariant under the slash-operator $\mid_{\kappa,\rho_L}(\gamma,\phi)$ for all $(\gamma,\phi) \in \widetilde{\Gamma}$, and if there exists a $\C[L'/L]$-valued Fourier polynomial (the \emph{principal part} of $f$)
	\begin{equation*}
	P_f(\tau) \coloneqq \sum_{r (2N)} \sum_{D \leq 0} c_f^+(D,r) e\left(\frac{D\tau}{4N}\right) \e_r 
	\end{equation*}
	such that $f(\tau) - P_f(\tau) = O(e^{-\varepsilon v})$ as $v \rightarrow \infty$ for some $\varepsilon > 0$. We denote the vector space of harmonic Maass forms of weight $\kappa$ with respect to $\rho_L$ by $H_{\kappa,\rho_L}$, and we let $M_{\kappa,\rho_L}^!$ be the subspace of weakly holomorphic modular forms. Every $f \in H_{\kappa,\rho_L}$ can be written as a sum $f = f^+ + f^-$ of a \emph{holomorphic} and a \emph{non-holomorphic part}, having Fourier expansions of the form
\begin{align*}
f^+(\tau) &= \sum_{r (2N)} \sum_{D \gg -\infty} c_f^+(D,r) e\left(\frac{D\tau}{4N}\right) \e_r, \\
f^-(\tau) &= \sum_{r (2N)} \sum_{D < 0} c_f^-(D,r)\Gamma\left(1-\kappa, \frac{\pi |D| v}{N}\right) e\left(\frac{D\tau}{4N}\right) \e_r,
\end{align*}
where $\Gamma(s,x): = \int_x^\infty t^{s-1}e^{-t} dt$ denotes the {\it incomplete Gamma function}.

The antilinear differential operator $\xi_\kappa := 2iv^\kappa \overline{\frac{\partial}{\partial {\overline{\tau}}}}$ maps a harmonic Maass form $f \in H_{\kappa,\rho_L}$ to a cusp form of weight  $2-\kappa$ for $\overline{\rho}_{L}$. We further require the \emph{raising operator} $R_\kappa := 2i\frac{\partial}{\partial \tau}+ \frac{\kappa}{v}$, which raises the weight of a smooth function transforming like a modular form of weight $\kappa$ for $\rho_L$ by two, and we define the \emph{iterated raising operator} $R_{\kappa}^n := R_{\kappa+2n-2}\circ \dots \circ R_{\kappa + 2}\circ R_\kappa$ for $n \in \N$, and $R_\kappa^0 := \mathrm{id}$.

\subsection{Maass Poincar\'e series}
\label{poincare series}

Examples of harmonic Maass forms can be constructed using Maass Poincar\'e series. We recall their construction from \cite[Section~1.3]{brhabil}. Let $\kappa \in \Z + \frac{1}{2}$ with $\kappa < 0$. For $s \in \C$ and $v > 0$ we let
\[
\mathcal{M}_{s}(v) := v^{-\kappa/2}M_{-\kappa/2,s-1/2}(v),
\]
with the usual $M$-Whittaker function. Let $D \in \Z$ with $D < 0$ and $r \in \Z/2N\Z$ with $D \equiv r^2 \pmod {4N}$. For $s \in \C$ with $\Re(s) > 1$ we define the $\C[L'/L]$-valued Maass Poincar\'e series
\[
\mathcal{P}_{\kappa,D,r}(\tau,s) := \frac{1}{2\Gamma(2s)}\sum_{(M,\phi) \in \widetilde{\Gamma}_{\infty}\backslash \Mp_{2}(\Z)}\mathcal{M}_{s}\left(\frac{\pi |D| v}{N}\right)e\left(\frac{Du}{4N}\right)\e_{r}\bigg|_{\kappa,\rho_L}(M,\phi),
\]
where $\widetilde{\Gamma}_{\infty}$ is the subgroup of $\Mp_{2}(\Z)$ generated by $T=\left(\left(\begin{smallmatrix}1 & 1 \\ 0 & 1 \end{smallmatrix} \right), 1 \right)$. It converges absolutely for $\Re(s) > 1$, it transforms like a modular form of weight $\kappa$ for $\rho_{L}$, and it is an eigenform of the Laplace operator $\Delta_{k}$ with eigenvalue $s(1-s)+(\kappa^{2}-2\kappa)/4$. The special value
\[
\mathcal{P}_{\kappa,D,r}(\tau) := \mathcal{P}_{\kappa,D,r}\left(\tau,1-\frac{\kappa}{2}\right)
\]
defines a harmonic Maass form in $H_{\kappa,\rho_L}$ whose Fourier expansion starts with
\[
\mathcal{P}_{\kappa,D,r}(\tau) = q^{-\frac{|D|}{4N}}(\e_{r}+(-1)^{\kappa - \frac{1}{2}}\e_{-r}) + \mathcal{O}(1).
\]
The following lemma follows inductively from \cite[Proposition~2.2]{bruinieronoalgebraic}.

\begin{lemma}\label{raising poincare series}
	For $n \in \N_0$ we have that
	\begin{equation*}
	R_\kappa^n \mathcal{P}_{\kappa,D,r}(\tau,s) = \left(\frac{\pi |D|}{N}\right)^n \frac{\Gamma\left(s+n+\frac{\kappa}{2}\right)}{\Gamma\left(s+\frac{\kappa}{2}\right)} \mathcal{P}_{\kappa+2n,D, r}(\tau,s).
	\end{equation*}
\end{lemma}

Maass Poincar\'e series for the dual Weil representation $\overline{\rho}_{L}$ are defined analogously, and Lemma~\ref{raising poincare series} remains true for them. However, in this case, we have to require that $D \equiv -r^2 \pmod{4N}$, and the Fourier expansion starts with $q^{-|D|/4N}(\e_r - (-1)^{\kappa-\frac{1}{2}}\e_{-r}) + \mathcal{O}(1)$.

\subsection{$W$-Whittaker functions}
We record some facts about derivatives and integrals of $W$-Whittaker functions that will be used in the computations later on. Following \cite[Section 1.3]{brhabil} we put
\[
\mathcal{W}_{\kappa,s}(y) := |y|^{-\kappa/2}W_{\frac{\kappa}{2}\sgn(y),s-\frac{1}{2}}(|y|), \qquad (\kappa \in \R, \ s \in \C, \ y \in \R \setminus \{0\}),
\]
where $W_{\nu,\mu}(y)$ is the usual $W$-Whittaker function. At $s = 1-\frac{\kappa}{2}$ it simplifies to
\begin{align}\label{special values W}
\mathcal{W}_{\kappa,1-\frac{\kappa}{2}}(y) = \begin{cases}
e^{-y/2}, & \text{if $y > 0$}, \\
e^{-y/2}\Gamma(1-\kappa,|y|), & \text{if $y < 0$}.
\end{cases}
\end{align}
Moreover, using (13.4.33) and (13.4.31) in \cite{abramowitz}, we obtain the formula
\begin{align}\label{raising W}
\begin{split}
&R_{\kappa}\left(\mathcal{W}_{\kappa,s}(4\pi m v)e(mu)\right) \\
&\qquad = \begin{cases}
-4\pi |m|\left(s+\frac{\kappa}{2}\right)\left(s-\frac{\kappa}{2}-1\right)\mathcal{W}_{\kappa+2,s}(4\pi m v)e(mu), & \text{if $m < 0$}, \\
-4\pi |m|\mathcal{W}_{\kappa+2,s}(4\pi m v)e(mu), & \text{if $m > 0$},
\end{cases}
\end{split}
\end{align}
for $u+iv \in \H$ and $m \in \R \setminus \{0\}$. Moreover, we will need the integral formula
\begin{align}\label{integral W}
\int_0^\infty v^{\kappa-2}\mathcal{W}_{\kappa,s}(\alpha v)\exp\left(-\frac{\alpha v}{2}-\frac{\beta}{v} \right)dv = \alpha^{\frac{1}{4}-\frac{\kappa}{2}}\beta^{\frac{\kappa}{2}-\frac{3}{4}}\sqrt{\pi}\mathcal{W}_{0,\frac{3}{2}-2s}\left( 4\sqrt{\alpha\beta}\right)
\end{align}
for $\alpha,\beta > 0$, which follows by evaluating the integral in terms of the $K$-Bessel function using \cite[(22) on p. 217]{tables}, and then writing $\sqrt{2z/\pi}K_{\mu}(z) = W_{0,\mu}(2z)$, see \cite[(13.6.21)]{abramowitz}.

\subsection{Hypergeometric series}
We compute the action of the iterated raising operator on a certain hypergeometric series for later use.
\begin{proposition}\label{raising hypergeometric}
For $X \in V(\R)$ with $q(X) = m < 0$ and $k \in \N$ we have
\begin{align*}
R_{0,z}^k &\bigg( \left(\frac{2|m|}{p_z^2(X)}\right)^{\frac{k}{2}}{}_2F_1\left( \frac{k}{2}, \frac{k+1}{2};k+\frac12;\frac{2|m|}{p_z^2(X)}\right)\bigg) = \frac{\Gamma(2k)}{\Gamma(k)}\left(\sqrt{4N|m|}\frac{\mathrm{sgn}(p_z(X))}{Q_X(z)}\right)^{k}.
\end{align*}
\end{proposition}

\begin{proof}
For brevity, we put $w := w(z) := \frac{2|m|}{p_z^2(X)}$. We prove the proposition by induction. To this end, for $j \in \Z$ and fixed $k\in \N$ we define the function
\[
f_j(z) := \left(y^{-2}\overline{Q_X(z)} \right)^j w^{\frac{k+j}{2}} {}_2F_1\left( \frac{k+j}{2}, \frac{k+j+1}{2};k+\frac12;w\right).
\]
Note that we want to compute $R_0^k f_0(z)$. We claim that
\begin{align}\label{raising fj}
R_{2j}f_j(z) = \frac{\sgn(p_z(X))}{\sqrt{4N|m|}}(k+j)f_{j+1}(z).
\end{align}
To prove this, we first use that for every holomorphic function $g:\H \to \C$, every smooth function $h:\H \to \C$, and every $\kappa,\ell \in \R$ we have the simple relation 
\[
R_{\ell-\kappa}\left(y^{\kappa}\overline{g(z)}\cdot h(z)\right) = y^{\kappa}\overline{g(z)}\cdot R_{\ell}h(z).
\]
Applying this with $\kappa = -2j, \ell = 0,$ and $g(z) = Q_X^j(z)$, we obtain
\[
R_{2j}f_j(z) = \left(y^{-2}\overline{Q_X(z)} \right)^j R_0\left(w^{\frac{k+j}{2}} {}_2F_1\left( \frac{k+j}{2}, \frac{k+j+1}{2};k+\frac12;w\right)\right).
\]
Next, we write $R_0 = 2i\frac{\partial}{\partial z}$ and use the formula $
\frac{\partial}{\partial z}\big(z^a {}_{2}F_1(a,b;c;z)\big) = a z^{a-1} {}_{2}F_1(a+1,b;c;z)$
(see (15.5.3) of NIST) to obtain
\begin{align*}
&R_0\left(w^{\frac{k+j}{2}} {}_2F_1\left( \frac{k+j}{2}, \frac{k+j+1}{2};k+\frac12;w\right)\right) \\
&\quad= i(k+j) w^{\frac{k+j}{2}-1} {}_2F_1\left( \frac{k+j}{2}+1, \frac{k+j+1}{2};k+\frac12;w\right)\cdot \frac{\partial}{\partial z}w.
\end{align*}
A direct computation shows that
\[
\frac{\partial}{\partial z} w = -\frac{i }{\sqrt{4N|m|}}\sgn(p_z(X))y^{-2}\overline{Q_X(z)}w^{\frac{3}{2}}.
\]
Finally, using ${}_2F_1(a,b;c;z) = {}_2F_1(b,a;c;z)$ and taking everything together, we obtain \eqref{raising fj}.

Note that the left-hand side in the proposition equals $R_0^k f_0(z)$. Applying \eqref{raising fj} inductively, we find that $R_0^k f_0(z)$ is a multiple of $f_k(z)$. More explicitly, we get 
\begin{align*}
R_0^k f_0(z) 
&= \frac{\sgn(p_z(X))^k}{\sqrt{4N|m|}^k}\frac{(2k-1)!}{(k-1)!}\left(y^{-2}\overline{Q_X(z)} \right)^k w^{k} {}_2F_1\left( k, k+\frac{1}{2};k+\frac12;w\right).
\end{align*}
We have ${}_2F_1(a,b;b;z)=(1-z)^{-a}$ which implies 
\begin{align*}
{}_2F_1\left(k,k+\frac12;k+\frac12;w\right) &= \left( 1-\frac{2|m|}{p_z^2(X)}\right)^{-k} =  p_z^2(X)^k\left(\frac{|Q_X(z)|^2}{2Ny^2} \right)^{-k},
\end{align*}
where we used that $m = q(X) = q(X_z)+q(X_{z^\perp}) = \frac{1}{4Ny^2}|Q_X(z)|^2 - \frac{1}{2}p_z^2(X)$, compare~\eqref{pq identity}. Taking everything together, we obtain the formula in the proposition.

\end{proof}

\section{Canonical meromorphic modular forms of weight $2k$}
\label{meromorphic modular forms}

In this section we prove Proposition~\ref{canonical modular form} for $\Gamma:=\Gamma_0(N)$ in a series of lemmas. Let $k \in \Z$ with $k \geq 2$. For completeness, we first show that the pairing in \eqref{pairing} is well-defined.

\begin{lemma}
\label{pairing path independent}
	For $f \in \mathcal{D}_{2k,\R}(\Gamma)$ and any closed geodesic $C_\gamma$ on $\Gamma \backslash \H$ the value of the pairing $(f,C_\gamma)$ is independent of the choice of the path.
\end{lemma}

\begin{proof}
	We have to show that, for any pole $\varrho \in \H$ of $f$ and any small closed loop $\ell_\varrho$ around $\varrho$, we have
	\begin{align}\label{residue vanishes}
	\Re\left( i^{k-1}\int_{\ell_\varrho}f(z)Q_\gamma(z,1)^{k-1}\right) = 0.
	\end{align}
	This essentially follows from the residue theorem. For computational convenience, we will use the fact that every function $F$ which is meromorphic near $\varrho$ has an \emph{elliptic expansion} in weight $\kappa \in \Z$ of the shape
	\[
	F(z) = (z -\overline{\varrho})^{-\kappa}\sum_{n \gg -\infty}c_{F,\varrho}(n)X_{\varrho}(z)^{n}, \qquad X_{\varrho}(z):= \frac{z-\varrho}{z-\overline{\varrho}},
	\] 
	with coefficients $c_{F,\varrho}(n) \in \C$, see Proposition~17 in Zagier's part of \cite{zagier123}.
	Here $\kappa \in \Z$ can be chosen freely, but the coefficients $c_{F,\varrho}(n) \in \C$ will depend on $\kappa$. In particular, we do not need to require that $F$ transforms like a modular form of weight $\kappa$. 
	
	If $G$ is another meromorphic function near $\varrho$ with an elliptic expansion in weight $2-\kappa$ and coefficients $c_{G,\varrho}(n)$, then it follows from the residue theorem (see \cite[Lemma 4.1]{anbs}) that
	\begin{align}\label{residue evaluation}
	\int_{\ell_\varrho}F(z)G(z)dz = \frac{\pi}{\Im(\varrho)}\sum_{n \in \Z}c_{F,\varrho}(n)c_{G,\varrho}(-n-1).
	\end{align}
	
	By definition of the space $\mathcal{D}_{2k,\R}(\Gamma)$, the elliptic expansion (in weight $2k$) of $f \in \mathcal{D}_{2k,\R}(\Gamma)$ near $\varrho$ is of the shape
	\begin{align}\label{f elliptic expansion}
	f(z) = a_{f,\varrho}(2i\Im(\varrho))^k(z-\overline{\varrho})^{-2k}X_{\varrho}(z)^{-k} + \mathcal{O}(1),
	\end{align}
	as $z \to \varrho$, with $a_{f,\varrho} \in \R$. Since $Q_\gamma(z,1)^{k-1}$ is holomorphic at $\varrho$, by \eqref{residue evaluation} we only need its elliptic expansion coefficient (in weight $2-2k$) of index $k-1$ at $\varrho$. By \cite[Lemma 3.1, Lemma 5.4]{anbs} this coefficient is given by
	\begin{align}
	\begin{split}\label{Q elliptic expansion}
	&\frac{(-4\Im(\varrho))^{1-k}}{(k-1)!}R_{2-2k}^{k-1}\big(Q_\gamma(z,1)^{k-1}\big)|_{z = \varrho}  \\
	&\quad= (-4\Im(\varrho))^{1-k}\left(2i\sqrt{d_\gamma}\right)^{k-1}P_{k-1}\left(\frac{i(A|\varrho|^2 + B\Re(\varrho) + C)}{\Im(\varrho)\sqrt{d_\gamma}} \right),
	\end{split}
	\end{align}
	where $d_\gamma$ is the discriminant of $Q_\gamma$, $P_{\ell}(x)$ denotes the $\ell$-th Legendre polynomial, and we wrote $Q_\gamma = [A,B,C]$. Using that $P_{\ell}$ is even if $\ell$ is even and odd if $\ell$ is odd, we see that \eqref{Q elliptic expansion} is real. Hence, putting \eqref{f elliptic expansion} and \eqref{Q elliptic expansion} in \eqref{residue evaluation}, we see that \eqref{residue evaluation} lies in $i^{k}\R$. This shows \eqref{residue vanishes}, and finishes the proof.
\end{proof}

For $\varrho \in \H$ and $k \in \Z$ with $k\geq 2$ we consider the \emph{Petersson Poincar\'e series} \cite{petersson}
\begin{align}\label{petersson poincare series}
\eta_{k,\varrho}(z) := \frac{1}{2}\sum_{\gamma \in \Gamma}\left( \frac{(z-\varrho)(z-\overline{\varrho})}{\varrho-\overline{\varrho}}\right)^{-k} \bigg|_{2k}\gamma.
\end{align}
It is straightforward to check that $\eta_{k,\varrho}(z)$ only depends on the class of $\varrho$ in $\Gamma \backslash \H$ and defines a meromorphic modular form in $\mathcal{D}_{2k,\R}(\Gamma)$ with residue divisor $\mathrm{res}(\eta_{k,\varrho}) = [\varrho]$. For a given divisor $D = \sum_{[\varrho]}\frac{c_{\varrho}}{w_\varrho}[\varrho] \in \mathrm{Div}(\Gamma \backslash \H)_{\R}$ we put
\[
\eta_{k,D}(z) := \sum_{[\varrho] \in \Gamma \backslash \H}\frac{c_{\varrho}}{w_\varrho}\eta_{k,\varrho}(z) \in \mathcal{D}_{2k,\R}(\Gamma).
\]
Then $\mathrm{res}(\eta_{k,D}) = D$. In order to prove Proposition~\ref{canonical modular form}, we need to show that $(\eta_{k,D},C) = 0$ for any $C \in \mathcal{C}_{\R}(\Gamma)$, and that $\eta_{k,D}$ is uniquely determined by this condition among all forms in $\mathcal{D}_{2k,\R}(\Gamma)$ with residue divisor $D$.

\begin{lemma}
	For any $C \in \mathcal{C}_{\R}(\Gamma)$ we have $(\eta_{k,D},C) = 0$.
\end{lemma}

\begin{proof}
	It suffices to prove the lemma for $D = \frac{1}{w_\varrho}[\varrho]$ and $C = C_\gamma$ being a closed geodesic on $\Gamma\backslash \H$. In this case the claim follows almost verbatim from the proof of \cite[Theorem~3]{katok}: Using the unfolding argument one obtains that $(\eta_{k,\varrho},C_\gamma)$ is a finite linear combination of integrals of the shape $\int_{-i\infty}^{i\infty}\frac{z^{k-1}dz}{(Az^2 + Bz + C)^k}$ for certain real binary quadratic forms $[A,B,C]$ of negative discriminant. Since $AC > 0$ for any such quadratic form, \cite[Lemma~2]{katok} shows that these integrals all vanish.
\end{proof}

\begin{lemma}
	If $f \in \mathcal{D}_{2k,\R}(\Gamma)$ satisfies $\mathrm{res}(f) = D$ and $(f,C) = 0$ for any $C \in \mathcal{C}_{\R}(\Gamma)$, then $f = \eta_{k,D}$.
\end{lemma}

\begin{proof}
	Under the above assumptions, $g:= f- \eta_{k,D}$ is a cusp form of weight $2k$ for $\Gamma$ with 
	\[
	(g,C_\gamma) = \Re\left(\left(\frac{i}{\sqrt{d_\gamma}}\right)^{k-1}\int_{C_\gamma}g(z)Q_\gamma(z,1)^{k-1}dz \right) = 0
	\]
	for any closed geodesic $C_\gamma$ on $\Gamma \backslash \H$. In other words, either all geodesic cycle integrals of $g$ are real, or they are all purely imaginary (depending on the parity of $k$). By \cite[Theorem~2]{katok} this implies $g = 0$. Alternatively, one may use the Eichler-Shimura isomorphism to deduce that $g = 0$.
\end{proof}

Taking together the above lemmas, we obtain the statement of Proposition~\ref{canonical modular form}. We remark that Proposition~\ref{canonical modular form} holds more generally for discrete subgroups $\Gamma \subset \SL_2(\R)$ with finite covolume, by the same arguments as above.

We close this section with some remarks about further properties of the meromorphic modular forms $\eta_\varrho(z)$.

\begin{remark}
	\begin{enumerate}
		\item One can show that $\eta_{k,\varrho}(z)$ is a constant multiple of $R_{0,z}^{k}G_{k}(z,\varrho)$, where $G_{k}(z,\varrho)$ denotes the higher Green function on $\H \times\H$, see \cite{mellit, bringmannkanevonpippich}.
		\item It is well-known that $\eta_{k,\varrho}(z)$ is orthogonal to cusp forms with respect to a suitable regularized Petersson inner product, see \cite{petersson, bringmannkanevonpippich}.
		\item We have seen above that the \emph{real parts} of the cycle integrals of $i^{k-1}\eta_{k,\varrho}(z)$ vanish. In \cite{anbs,loebrichschwagenscheidt} it was proved that if $\varrho$ is a CM point, then the \emph{imaginary parts} of (certain linear combinations of) the cycle integrals of $i^{k-1}\eta_{k,\varrho}(z)$ are \emph{rational} (up to some normalizing factors involving powers of $\pi$ and the disciminant of $\varrho$). 
	\end{enumerate}
\end{remark}

\section{Theta lifts and meromorphic modular forms}
\label{theta lifts}

\subsection{Twisted theta functions and theta lifts} Let $L$ be the lattice from Section~\ref{quadratic spaces}. Let $\Delta \in \Z$ be a fundamental discriminant (possibly $1$) and let $\rho \in \Z$ be such that $\rho^2 \equiv \Delta \pmod{2N}$. We let $\chi_\Delta$ be the generalized genus character on binary quadratic forms as in \cite{gkz}, which can also be viewed as a function on $L'$ using the identification of $L_{d,r}$ with $\mathcal{Q}_{d,r}$ explained in Section~\ref{quadratic spaces}. For notational convenience we put
\[
\widetilde{\rho}_L := \begin{cases}
\rho_L, & \text{if } \Delta > 0, \\
\overline{\rho}_{L}, & \text{if } \Delta < 0.
\end{cases}
\]

Following \cite{bruinieronoalgebraic}, we define the \emph{twisted Siegel theta function} by
\[
\Theta_{\Delta,\rho}(\tau,z) := v^\frac{1}{2}\sum_{h \in L'/L}\sum_{\substack{X \in L+\rho h \\ q(X) \equiv \Delta q(h)\, (\Delta)}}\chi_\Delta(X)e\left(\frac{q(X_z)}{|\Delta|}\tau + \frac{q(X_{z^\perp})}{|\Delta|}\overline{\tau} \right)\e_h.
\]
It transforms like a modular form of weight $\frac{1}{2}$ for $\widetilde{\rho}_L$ in $\tau$ and is $\Gamma_0(N)$-invariant in $z$. Similarly, following \cite{hoevel} we define the \emph{twisted Millson theta function} by
\[
\Theta^*_{\Delta,\rho}(\tau,z) := v^\frac{3}{2}\sum_{h \in L'/L}\sum_{\substack{X \in L+\rho h \\ q(X) \equiv \Delta q(h)\, (\Delta)}}\chi_\Delta(X)p_z(X)e\left(\frac{q(X_z)}{|\Delta|}\tau + \frac{q(X_{z^\perp})}{|\Delta|}\overline{\tau} \right)\e_h,
\]
with the polynomial $p_z(X) = -(X,X(z))$ as in \eqref{polynomials}. The Millson theta function transforms like a modular form of weight $-\frac{1}{2}$ for $\widetilde{\rho}_L$ in $\tau$ and is $\Gamma_0(N)$-invariant in $z$. These theta functions and the corrsponding theta lifts have been studied at many places in the past years, see for example \cite{asmillson, bruinierehlenyang, crawfordfunke, hoevel}.

Let $k \in \N$ with $k \geq 2$. For a harmonic Maass form $f \in H_{\frac{3}{2}-k,\widetilde{\rho}_L}$ we define the theta lift
\begin{align*}
\Phi_{\Delta,\rho}^k(f,z) := \begin{dcases}
C_k^- \cdot R_{0,z}^{k}\int_\mathcal{F}^\reg \left\langle R_{\frac{3}{2}-k,\tau}^{\frac{k-1}{2}}f(\tau),\Theta_{\Delta,\rho}(\tau,z)\right\rangle v^\frac{1}{2} d\mu(\tau), & \text{if $k$ is odd}, \\
C_k^+ \cdot R_{0,z}^{k}\int_\mathcal{F}^\reg \left\langle R_{\frac{3}{2}-k,\tau}^{\frac{k}{2}-1}f(\tau),\Theta_{\Delta,\rho}^*(\tau,z)\right\rangle v^{-\frac{1}{2}} d\mu(\tau), & \text{if $k$ is even},
\end{dcases}
\end{align*}
where $d\mu(\tau) := \frac{du dv}{v^2}$ denotes the usual invariant measure on $\H$, the integrals are regularized as in \cite[Section~6]{borcherds}, and $C_k^-$ and $C_k^+$ denote the normalizing constants
\begin{align}\label{normalizing constants}
C_{k}^- := \frac{i^k|\Delta|^{\frac{k-1}{2}}N^{\frac{k-1}{2}}}{2^{k+1}\pi^{\frac{k-1}{2}}  (k-1)!}, \qquad C_k^+ := \frac{i^k|\Delta|^{\frac{k}{2}-1}N^{\frac{k-1}{2}}}{2^{k+\frac{1}{2}}\pi^{\frac{k}{2}-1} (k-1)!}.
\end{align}
We remark that Zemel \cite{zemel} investigated a similar lift for lattices of more general signature (with $\Delta = 1$) and used it to prove a higher weight version of the Gross-Kohnen-Zagier Theorem.  Moreover, an analogous theta lift (defined in the same way as $\Phi_{\Delta,\rho}^k(f,z)$ but without the raising operator $R_{0,z}^k$ in front) was considered in \cite{bruinierehlenyang}, where it was shown to be a linear combination of higher Green functions $G_k(Z_{\Delta,\rho}(D,r),z)$ evaluated at Heegner divisors $Z_{\Delta,\rho}(D,r)$ in the first variable.

Here, we will show that the lift $\Phi_{\Delta,\rho}^k(f,z)$ defines a meromorphic modular form of weight $2k$ for $\Gamma_0(N)$ with poles of order $k$ at a linear combination of Heegner divisors $Z_{\Delta,\rho}(D,r)$. Moreover, we will compute the Fourier expansion of the lift to study the algebraicity of the Fourier coefficients of these meromorphic modular forms.

\subsection{The theta lift as a meromorphic modular form}

Let $D < 0$ be an integer and let $r \in \Z/2N\Z$ with $D \equiv \sgn(\Delta)r^2 \pmod{4N}$. We compute the theta lift of the Maass Poincar\'{e} series $\mathcal{P}_{\frac{3}{2}-k,D,r}(\tau)$, defined in Section~\ref{poincare series}, and obtain an explicit representation of the lift in terms of the meromorphic modular form
\begin{align}\label{fkD}
f_{k,D,r,\Delta,\rho}(z) := i^k \big(|D\Delta|\big)^{k-\frac{1}{2}}\sum_{Q \in \mathcal{Q}_{D|\Delta|,r\rho}}\frac{\sgn(Q)\chi_{\Delta}(Q)}{Q(z,1)^{k}}
\end{align}
of weight $2k$ for $\Gamma_0(N)$. Here we put $\sgn(Q) := \sgn(a)$ for $Q = [aN,b,c] \in \mathcal{Q}_{D|\Delta|,r\rho}$. Let $\alpha_Q \in \H$ be the Heegner point corresponding to $Q \in \mathcal{Q}_{D|\Delta|,r\rho}$. Then a short computation using 
\[
Q(z,1) = 2i\sqrt{|D\Delta|}\frac{(z-\alpha_Q)(z-\overline{\alpha}_Q)}{(\alpha_Q-\overline{\alpha}_Q)}
\]
shows that
\[
f_{k,D,r,\Delta,\rho}(z) = (|D\Delta|)^{\frac{k-1}{2}}\sum_{Q \in \mathcal{Q}_{D|\Delta|,r\rho}/\Gamma}\sgn(Q)\chi_\Delta(Q)\eta_{k,\alpha_Q}(z),
\]
with the Petersson Poincar\'e series $\eta_{k,\varrho}(z)$ defined in \eqref{petersson poincare series}. In particular, $f_{k,D,r,\Delta,\rho} \in \mathcal{D}_{2k,\R}(\Gamma)$ is the canonical meromorphic modular form of weight $2k$ for the rescaled Heegner divisor 
\[
\mathrm{res}(f_{k,D,r,\Delta,\rho}) = |D\Delta|^{\frac{k-1}{2}}Z_{\Delta,\rho}(D,r).
\]

\begin{proposition}\label{theta unfolding}
Let $k \in \Z$ with $k \geq 2$, let $D \in \Z$ with $D < 0$ and $r \in \Z/2N\Z$ with $D \equiv \sgn(\Delta)r^2 \pmod{4N}$. Then we have
\begin{align*}
\Phi_{\Delta,\rho}^{k}\left(\mathcal{P}_{\frac{3}{2}-k,D,r}(\tau),z\right) &= 
f_{k,D,r,\Delta,\rho}(z).
\end{align*}
\end{proposition}

\begin{remark}
Let $f$ be a harmonic weak Maass form of weight $\frac{3}{2}-k$.
Since the Poincar\'e series $\mathcal{P}_{\frac{3}{2}-k,D,r}(\tau)$ generate the space of weight $\frac{3}{2}-k$ harmonic weak Maass forms, we directly deduce that the lift  $\Phi_{\Delta,\rho}^k(f,z)$ defines a meromorphic modular form of weight $2k$ for $\Gamma_0(N)$ with poles of order $k$ at a linear combination of Heegner divisors $Z_{\Delta,\rho}(D,r)$. 
\end{remark}

\begin{proof}
Let $k \geq 3$ be odd. Using Lemma~\ref{raising poincare series} we can write
\begin{align*}
&\Phi_{\Delta,\rho}^{k}\left(\mathcal{P}_{\frac{3}{2}-k,D,r}(\tau),z\right) \\
&\quad = C_k^-\left( \frac{\pi|D|}{N}\right)^{\frac{k-1}{2}}\Gamma\left( \frac{k+1}{2}\right) R_{0,z}^{k}\int_\mathcal{F}^\reg \left\langle \mathcal{P}_{\frac{1}{2},D,r}\left(\tau,\frac{k}{2}+\frac{1}{4}\right),\Theta_{\Delta,\rho}(\tau,z)\right\rangle v^\frac{1}{2} d\mu(\tau),
\end{align*}
where $C_k^-$ is the normalizing constant defined in \eqref{normalizing constants}. The regularized integral can be computed by unfolding against the Maass Poincar\'e series as in the proof of \cite[Theorem~2.14]{brhabil}. We obtain that
\begin{align*}
&\Phi_{\Delta,\rho}^{k}\left(\mathcal{P}_{\frac{3}{2}-k,D,r}(\tau),z\right) = 2 C_k^- \left(\frac{\pi|D|}{N}\right)^{\frac{k-1}{2}} \frac{\Gamma\left(\frac{k+1}{2}\right)\Gamma\left(\frac{k}{2}\right)}{\Gamma\left(k+\frac12\right)}\\
& \quad \times \!\!\!\! \sum_{\substack{X \in L+\rho r\\q(X)=D|\Delta|/4N}} \chi_\Delta(X) \cdot R_{0,z}^k \left(    \left(\frac{D|\Delta|}{4Nq(X_{z^\perp})}\right)^{\frac{k}{2}}{}_2F_1\left( \frac{k}{2}, \frac{k+1}{2};k+\frac12;\frac{D|\Delta|}{4Nq(X_{z^\perp})}\right)\right).
\end{align*}
Recall from \eqref{pq identity} that $q(X_{z^\perp}) = -\frac{1}{2}p_z^2(X)$. By Proposition \ref{raising hypergeometric} (with $m = \frac{D|\Delta|}{4N}$) we find
\begin{align*}
\Phi_{\Delta,\rho}^{k}\left(\mathcal{P}_{\frac{3}{2}-k,D,r}(\tau),z\right) &= 2C_k^-\left(|D||\Delta|\right)^{\frac{k}{2}}\left(\frac{\pi|D|}{N}\right)^{\frac{k-1}{2}}\frac{\Gamma(2k)\Gamma\left(\frac{k+1}{2}\right)\Gamma\left(\frac{k}{2}\right)}{\Gamma(k)\Gamma\left(k+\frac12\right)} \\
 &\qquad  \times \!\!\!\! \sum_{\substack{X \in L+\rho r\\q(X)= D|\Delta|/4N}}  \mathrm{sgn}(p_z(X)) \frac{ \chi_\Delta(X)}{Q_X(z)^k}.
\end{align*}
The quotient of Gamma functions can be simplified to $2^k (k-1)!$ using the Legendre duplication formula $\sqrt{\pi}\Gamma(2z) = 2^{2z-1}\Gamma(z)\Gamma(z+\frac{1}{2})$. Moreover, we have $\sgn(p_z(X)) = \sgn(Q_X)$, independently of $z \in \H$, where $Q_X$ is the quadratic form corresponding to $X$ under the identification of $L_{D|\Delta|,r\rho}$ with $\mathcal{Q}_{N,D|\Delta|,r\rho}$ explained in Section~\ref{quadratic spaces}. Hence we can rewrite the sum into a multiple of $f_{k,D,r,\Delta,\rho}(z)$. Taking into account the normalizing constant $C_k^-$ defined in \eqref{normalizing constants} and putting everything together, we obtain the stated formula.

The proof for even $k$ is analogous, so we omit the details.
\end{proof}

\subsection{The Fourier expansion of the theta lift}
In this section we compute the Fourier expansion of the theta lift of harmonic Maass forms.
\begin{proposition}\label{theta lift fourier expansion}
	Let $k \in \Z$ with $k \geq 2$. For $y \gg 0$ large enough, the theta lift $\Phi_{\Delta,\rho}^k(f,z)$ of $f \in H_{\frac{3}{2}-k,\widetilde{\rho}_L}$ has the Fourier expansion
	\begin{align*}
	\Phi_{\Delta,\rho}^k(f,z) 
	= C_{k,\Delta} i \pi^k \sqrt{\Delta} \sum_{n \geq 1}n^{2k-1}\sum_{d \mid n}\left(\frac{\Delta}{d} \right) d^{-k}c_f^+\left(\frac{|\Delta|n^2}{d^2},\frac{\rho n}{d}\right) e(n\tau),
\end{align*}
with the rational constant $C_{k,\Delta} := \frac{(-2\sgn(\Delta))^k |\Delta|^{k-1}}{(k-1)!}$.
\end{proposition}

\begin{proof}
Suppose that $k \geq 3$ is odd. For simplicity we assume that $\Delta \neq 1$, but the proof for $\Delta = 1$ is very similar. We write the theta function $\Theta_{\Delta,\rho}(\tau,z)$ as a Poincar\'e series as in \cite[Theorem 4.8]{bruinieronoheegner}. Then we use the unfolding argument as in the proof of \cite[Theorem 5.3]{bruinieronoheegner} and obtain
\begin{align*}
 \Phi_{\Delta,\rho}^k(f,z)&=R_{0,z}^k\bigg(2\sqrt{N}C_k^-y\varepsilon\sum_{X \in K'}\sum_{n \geq 1}\left(\frac{\Delta}{n} \right)e(\sgn(\Delta)n(X,\mu)) \\
& \qquad \qquad   \times \int_0^\infty c(|\Delta|q(X),\rho X,v)\exp\left(-\frac{N\pi n^2y^2}{|\Delta|v}-2\pi q(X)|\Delta|v \right)\frac{dv}{v^{\frac{3}{2}}}\bigg),
\end{align*}
where $C_k^-$ is the normalizing constant defined in \eqref{normalizing constants}, $K' = \left\{\frac{m}{2N}\left(\begin{smallmatrix}1 & 0 \\ 0 & -1 \end{smallmatrix}\right) : m \in \Z\right\}$ is a positive definite sublattice of $L'$, $\varepsilon$ equals $1$ if $\Delta > 0$ and $i$ if $\Delta < 0$, we put $\mu=\left(\begin{smallmatrix} x & -x^2 \\ -1 & -x\end{smallmatrix}\right)$, and we wrote\footnote{note that, in contrast to \cite[Theorem~5.3]{bruinieronoheegner}, we put $e^{-2\pi n v}$ into $c(n,h,v)$.}
\[
R_{\frac{3}{2}-k}^{\frac{k-1}{2}}f(\tau) = \sum_{h \in L'/L}\sum_{n \gg -\infty}c(n,h,v)e(nu)\e_h.
\]
for the Fourier expansion of $R_{\frac{3}{2}-k}^{\frac{k-1}{2}}f(\tau)$.

We first show that the contribution for $X = 0$ vanishes. Note that 
\begin{align}\label{raising imaginary part}
R_\kappa v^\alpha = (\alpha+\kappa)v^{\alpha-1}
\end{align}
for any $\alpha,\kappa \in \R$, so we have
\[
c(0,0,v) = C_0c_f^+(0,0)v^{\frac{1-k}{2}}, \qquad \text{ where } C_0:= \prod_{j=1}^{\frac{k-1}{2}} \left(j-\frac{5}{2}\right).
\]
Now the contribution from the $X = 0$ summand equals
\begin{align*}
&2\sqrt{N}\varepsilon  C_k^-C_0c_f^+(0,0) \sum_{n \geq 1}\left(\frac{\Delta}{n} \right)R_{0,z}^k\left(y\int_0^\infty v^{\frac{1-k}{2}}\exp\left(-\frac{N\pi n^2y^2}{|\Delta|v}\right)\frac{dv}{v^{\frac{3}{2}}}\right) \\
&\quad= 2N^{\frac{1-k}{2}}\pi^{-\frac{k}{2}}|\Delta|^{\frac{k}{2}}\varepsilon  C_k^-C_0c_f^+(0,0)\Gamma\left(\tfrac{k}{2}\right)L_\Delta(k)\cdot R_{0,z}^k y^{1-k}.
\end{align*}
It follows from \eqref{raising imaginary part} that $R_{0}^k y^{1-k} = (k-1)! R_{2k-2}y^{2-2k} = 0$. In particular, the Fourier expansion of the theta lift has no constant term.

For $m > 0$ we have $e^{-2\pi m v} = \mathcal{W}_{\frac{3}{2}-k,\frac{k}{2}+\frac{1}{4}}(4\pi m v)$ by \eqref{special values W}, and hence
\begin{align*}
c(m,h,v)e(mu) &= c_f^+(m,h) R_{\frac{3}{2}-k}^{\frac{k-1}{2}}\left(\mathcal{W}_{\frac{3}{2}-k,\frac{k}{2}+\frac{1}{4}}(4\pi m v)e(mu) \right) \\
&=c_f^+(m,h)(-4\pi m)^{\frac{k-1}{2}}\mathcal{W}_{\frac{1}{2},\frac{k}{2}+\frac{1}{4}}(4\pi m v)e(mu)
\end{align*}
by \eqref{raising W}. Using \eqref{integral W} (with $\kappa = \frac{1}{2}, s = \frac{k}{2}+\frac{1}{4}, \alpha = 4\pi|\Delta|q(X), \beta = \frac{N\pi n^2y^2}{|\Delta|}$) we have for $X \neq 0$ the evaluation
\begin{align*}
&\int_0^\infty \mathcal{W}_{\frac{1}{2},\frac{k}{2}+\frac{1}{4}}(4\pi |\Delta|q(X) v)\exp\left(-\frac{N\pi n^2y^2}{|\Delta|v}-2\pi q(X)|\Delta|v \right)\frac{dv}{v^{\frac{3}{2}}}  \\
&\qquad = \frac{\sqrt{|\Delta|}}{\sqrt{N}ny}\mathcal{W}_{0,1-k}\left(4\pi ny \sqrt{4Nq(X)}\right).
\end{align*}
Note that $q(X) = \frac{m^2}{4N}$ and $(X,\mu) = mx$ for $X = \frac{m}{2N}\left(\begin{smallmatrix}1 & 0 \\ 0 & -1 \end{smallmatrix} \right)\in K'$. Using \eqref{raising W} we find 
\begin{align*}
&R_0^k\big(\mathcal{W}_{0,1-k}(4\pi |m|ny)e(\sgn(\Delta)mnx)\big)\\
&\qquad = \begin{cases} 
(-4\pi|m|n)^k\mathcal{W}_{k,1-k}(4\pi |m|ny)e(|m|nx), & \sgn(\Delta)m > 0, \\
0, & \sgn(\Delta)m < 0, 
\end{cases}
\end{align*}
since for $\sgn(\Delta)m < 0$ one of the factors $(s+\frac{\kappa}{2})$ in \eqref{raising W} will be $0$. By \eqref{special values W} we have $\mathcal{W}_{k,1-k}(4\pi |m|n y) = e^{-2\pi |m|ny}$. Taking everything together, we obtain after simplification
\begin{align*}
&\Phi_{\Delta,\rho}^k(f,z) \\
&= (-1)^{\frac{k+1}{2}}\varepsilon C_k^- \sgn(\Delta)2^{2k+1}N^{\frac{1-k}{2}}\pi^{\frac{3k-1}{2}}|\Delta|^{\frac{k}{2}}\sum_{n \geq 1}n^{2k-1}\sum_{d \mid n}\left(\frac{\Delta}{d} \right) d^{-k}c_f^+\left(\frac{|\Delta|n^2}{d^2},\frac{\rho n}{d}\right) e(n\tau),
\end{align*}
where we also used that $c_f^+(m,\sgn(\Delta)r) = \sgn(\Delta)c_f^+(m,r)$. Putting in the normalizing constant $C_k^-$ from \eqref{normalizing constants}, we obtain the stated formula.

The proof for even $k$ is similar. Now we use Theorem~5.10 and Lemma~5.6 of \cite{crawfordfunke} to write the theta function $\Theta_{\Delta,\rho}^*(\tau,z)$ as a Poincar\'e series
\begin{align*}
\Theta_{\Delta,\rho}^*(\tau,z)&=  -\frac{Niy^2}{2\sqrt{2|\Delta|}} \sum_{n\geq 1} \sum_{\tilde\gamma\in \tilde\Gamma_\infty \setminus\tilde\Gamma} n\left[ \exp\left( -\frac{N\pi n^2y^2}{|\Delta| v}\right) \Xi(\tau,\mu,-n,0)\right]_{-1/2,\tilde\rho_K}\tilde\gamma,
\end{align*}
where
\[
\Xi(\tau,\mu,-n,0) = \left(\frac{\Delta}{n}\right)\varepsilon\sqrt{|\Delta|} \sum_{h \in K'/K} \sum_{\substack{X \in K+rh\\Q(\lambda)\equiv \Delta Q(h)(\Delta)}} e\left( \frac{q(X) \tau}{|\Delta|}- \frac{n(X,\mu)}{|\Delta|}\right)\e_h.
\] 
By the unfolding argument we obtain as in the proof of \cite[Theorem~5.3]{bruinieronoheegner} that
\begin{align*}
 \Phi_{\Delta,\rho}^k(f,z)&=R_{0,z}^k\bigg(\sqrt{2}Niy^2\bar\varepsilon \sum_{X \in K'}\sum_{n \geq 1}n \left(\frac{\Delta}{n} \right) e(\sgn(\Delta)n(\lambda,\mu)) \\
& \qquad \qquad   \times \int_0^\infty c(|\Delta|q(X),\rho X,v)\exp\left(-\frac{N\pi n^2y^2}{|\Delta|v}-2\pi q(X)|\Delta|v \right)\frac{dv}{v^{\frac{5}{2}}}\bigg).
\end{align*}
Now the rest of the proof proceeds as in the case of odd $k$, so we omit the details.
\end{proof}

\section{Meromorphic modular forms corresponding to Heegner divisors}
\label{meromorphic modular forms for heegner divisors}

In this section we investigate the canonical and normalized meromorphic modular forms corresponding to twisted Heegner divisors and prove Theorem~\ref{canonical modular form heegner} and Theorem~\ref{normalized modular form heegner}.

Let $N \in \N$ and $k \in \N$ with $k \geq 2$, and let $L$ be the lattice defined in \eqref{the lattice}. Recall from \cite{ez,skoruppa} that the space $M_{\frac{1}{2}+k,\rho_L}$ of holomorphic vector-valued modular forms of weight $\frac{1}{2}+k$ for $\rho_L$ is isomorphic to the space $J_{k+1,N}^{\mathrm{skew}}$ of skew holomorphic Jacobi forms of weight $k+1$ and index $N$. Similarly, $M_{\frac{1}{2}+k,\overline{\rho}_L}$ is isomomorphic to the space $J_{k+1,N}$ of holomorphic Jacobi forms of weight $k+1$ and index $N$. Hence, the theory of Hecke operators for Jacobi forms developed in \cite{ez,skoruppazagier,skoruppa} carries over to vector-valued modular forms for $\rho_L$ and $\overline{\rho}_L$. In particular, there is a newform theory for these spaces of vector-valued modular forms.

We let $S_{2k}^{\new,+}(N)$ be the space of cuspidal newforms of weight $2k$ for $\Gamma_0(N)$ on which the Fricke involution acts by multiplication with $(-1)^k$. 

By the Shimura correspondence, it is isomorphic as a module over the Hecke algebra to  $S_{\frac{1}{2}+k,\rho_L}^{\new}$, compare \cite{skoruppazagier,gkz,skoruppa}. Similarly, the space $S_{2k}^{\new,-}(N)$ of cuspidal newforms of weight $2k$ for $\Gamma_0(N)$ on which the Fricke involution acts by multiplication with $-(-1)^k$ is isomorphic as a module over the Hecke algebra to $S_{\frac{1}{2}+k,\overline{\rho}_L}^{\new}$.

Let $\rho$ be one of the representations $\rho_L$ and $\overline{\rho}_L$. For every $m \in \N$ there is a Hecke operator $T_m$ acting on $M_{\frac{1}{2}+k,\rho}$. The action on Fourier expansions can be computed explicitly, compare \cite{ez,skoruppazagier}. For example, if $p$ is a prime with $(p,N) = 1$, and we let $f = \sum_{n,r}c_f(n,r)q^n\e_r \in M_{\frac{1}{2}+k,\rho}$ and $f|T_p = \sum_{n,r}c_{f|T_p}(n,r)q^n\e_r$, then we have
\begin{align}\label{hecke operator}
c_{f|T_p}(n,r) = c_{f}(p^2 n,pr) + p^{k-1}\left(\frac{4N\sigma n}{p} \right)c_f(n,r) + p^{2k-1}c_f(n/p^2,r/p),
\end{align}
where $\sigma = 1$ if $\rho = \rho_L$ and $\sigma = -1$ if $\rho = \overline{\rho}_L$. There are similar formulas for general $m \in \N$. The Hecke operators act on harmonic Maass forms in an analogous way, and the action on Fourier expansions is the same.

We let $\rho \in \Z/2N\Z$ and $\Delta \equiv \rho^2 \pmod{4N}$ be a fundamental discriminant, and we put $\epsilon := -\sgn(\Delta)$. For a normalized newform $G \in S_{2k}^{\new,\epsilon}(N)$ we let $F_G$ be the totally real number field generated by its Fourier coefficients, and we let $g \in S_{\frac{1}{2}+k,\overline{\widetilde{\rho}}_L}^{\new}$ be its Shimura correspondent, which we may normalize to have coefficients in $F_G$, as well. Then, by \cite[Lemma~7.2]{bruinieronoheegner} there exists a harmonic Maass form $f \in H_{\frac{3}{2}-k,\widetilde{\rho}_L}(F_G)$ whose principal part has coefficients in $F_G$, and which satisfies $\xi_{\frac{3}{2}-k}f = \|g\|^{-2}g$, where we write $\|g\| = \sqrt{(g,g)}$ for the Petersson norm of $g$. We let
\[
Z_{k,\Delta,\rho}(f) := \sum_{r(2N)}\sum_{D < 0}c_f^+(D,r)|D\Delta|^{\frac{k-1}{2}}Z_{\Delta,\rho}(D,r)
\]
be the twisted Heegner divisor corresponding to $f$. Then the canonical meromorphic modular form with residue divisor $Z_{k,\Delta,\rho}(f)$ is given by
\begin{align}\label{canonical modular form Zf}
\eta_{k,\Delta,\rho}(f) = \sum_{r(2N)}\sum_{D < 0}c_f^+(D,r)f_{k,D,r,\Delta,\rho}(z),
\end{align}
with $f_{k,D,r,\Delta,\rho}$ defined in \eqref{fkD}. On the other hand, by Proposition~\ref{theta unfolding} the meromorphic modular form $\eta_{k,\Delta,\rho}(f)$ can be obtained as the theta lift $\Phi_{\Delta,\rho}^k(f,z)$ of $f$. In particular, by Proposition~\ref{theta lift fourier expansion}, we have the Fourier expansion
\begin{align}\label{fourier expansion canonical modular form Zf}
\eta_{k,\Delta,\rho}(f) = C_{k,\Delta}i\pi^k\sqrt{\Delta}\sum_{n > 0}n^{2k-1}\sum_{d \mid n}\left( \frac{\Delta}{d}\right)d^{-k} c_f^+\left(\frac{|\Delta|n^2}{d^2}, \frac{\rho n}{d}\right)e(nz)
\end{align}
with a rational constant $C_{k,\Delta} \in \Q$.

\begin{theorem}\label{canonical modular form heegner level N}
	Let $f \in H_{\frac{3}{2}-k,\widetilde{\rho}_L}(F_G)$ be as above, and let $\eta_{k,\Delta,\rho}(f) \in \mathcal{D}_{2k,\R}(\Gamma)$ be the canonical meromorphic modular form of weight $2k$ for $Z_{k,\Delta,\rho}(f)$ as in \eqref{canonical modular form Zf}. Then the following are equivalent.
	\begin{enumerate}
		\item We have $c_{f}^+(|\Delta|,\rho) \in F_G$.
		\item We have $c_{f}^+(n^2|\Delta|,n\rho) \in F_G$ for all $n \in \N$.
		\item All Fourier coefficients of $\eta_{k,\Delta,\rho}(f)$ are contained in $i\pi^k \sqrt{\Delta} F_G$.
	\end{enumerate}
 \end{theorem}
 
 \begin{proof}
 	The implication $(2)$ $\Rightarrow$ $(3)$ follows from the Fourier expansion of $\eta_{k,\Delta,\rho}(f)$ in \eqref{fourier expansion canonical modular form Zf}. Similarly, the implication $(3)$ $\Rightarrow$ $(1)$ follows from \eqref{fourier expansion canonical modular form Zf} since the first coefficient of $\eta_{k,\Delta,\rho}$ is given by $C_{k,\Delta}i\pi^k\sqrt{\Delta}c_f^+(|\Delta|,\rho)$. 
	
	It remains to prove that $(1)$ implies $(2)$. We use the action of the Hecke operators on $f$. By \cite[Proposition~7.1]{bruinieronoheegner} we have $\xi_{\frac{3}{2}-k}(f|T_n) = n^{1-2k}\xi_{\frac{3}{2}-k}(f)|T_n$. Moreover, using $\xi_{\frac{3}{2}-k}(f) = \|g\|^{-2}g$ and $g|T_{n} = \lambda_n(G)g$, where $\lambda_n(G) \in F_G$ is the $T_n$-eigenvalue of $G$, we see that
\[
f' := n^{2k-1}f|T_n - \lambda_n(G)f \in M_{\frac{3}{2}-k,\widetilde{\rho}_L}^!(F_G)
\]
is a weakly holomorphic modular form with principal part defined over $F_G$. Since the space of weakly holomorphic modular forms for $\widetilde{\rho}_L$ has a basis of forms with rational Fourier coefficients by a result of McGraw \cite{mcgraw}, we obtain that all Fourier coefficients of $f'$ lie in $F_G$ . Now using the explicit formula \eqref{hecke operator} for the action of $T_n$ on the coefficients of $f$, we see that $c_f^+(n^2|\Delta|,n\rho)$ is a rational linear combination of coefficients $c_f^+(m^2 |\Delta|,m\rho)$ for some $m < n$ and coefficients of $f'$, which lie in $F_G$ by the above discussion. Hence, if $c_f^+(|\Delta|,\rho) \in F_G$ then it follows by induction that $c_f^+(n^2|\Delta|,n\rho) \in F_G$ for every $n \in \N$.
 \end{proof}

\begin{theorem}\label{normalized modular form heegner level N}
	For $f \in H_{\frac{3}{2}-k,\widetilde{\rho}_L}(F_G)$ as above, there exists a unique meromorphic modular form $\zeta_{k,\Delta,\rho}(f) \in \mathcal{D}_{2k,\R}(\Gamma)$ with the following properties.
	\begin{enumerate}
		\item The residue divisor is given by $\mathrm{res}(\zeta_{k,\Delta,\rho}(f)) = Z_{k,\Delta,\rho}(f)$.
		\item We have $(\zeta_{k,\Delta,\rho}(f),C) = 0$ for every closed geodesic $C$ in $\Gamma_0(N) \backslash \H$ with $(G,C) = 0$.			
		\item The first Fourier coefficient of $\zeta_{k,\Delta,\rho}(f)$ vanishes.
	\end{enumerate}
	Moreover, all Fourier coefficients of $\zeta_{k,\Delta,\rho}(f)$ lie in $i\pi^k \sqrt{\Delta}F_G$, and we have
	\begin{align}\label{coefficient period formula}
	c_f^+(|\Delta|,\rho) = -\frac{1}{C_{k,\Delta}\pi i^k\sqrt{\Delta}}\cdot\frac{(\zeta_{k,\Delta,\rho}(f),C)}{(G,C)}
	\end{align}
	for every $C \in \mathcal{C}_{\R}(\Gamma)$ with $(G,C) \neq 0$.
\end{theorem}

\begin{proof}
	By \eqref{fourier expansion canonical modular form Zf} the first coefficient of $\eta_{k,\Delta,\rho}(f)$ is given by $C_{k,\Delta}i\pi^k\sqrt{\Delta}c_f^+(|\Delta|,\rho)$. Hence
	\[
	\zeta_{k,\Delta,\rho}(f):= \eta_{k,\Delta,\rho}(f) - C_{k,\Delta}i\pi^k\sqrt{\Delta}c_f^+(|\Delta|,\rho)G
	\]
	has the desired properties. Taking the bilinear pairing with any $C \in \mathcal{C}_{\R}(\Gamma)$ yields
	\begin{align*}
	(\zeta_{k,\Delta,\rho}(f),C) &= (\eta_{k,\Delta,\rho}(f),C) - C_{k,\Delta}i\pi^k\sqrt{\Delta}c_f^+(|\Delta|,\rho)(G,C) \\
	&= -C_{k,\Delta}i\pi^k\sqrt{\Delta}c_f^+(|\Delta|,\rho)(G,C),
	\end{align*}
	which gives the formula \eqref{coefficient period formula}.
	
	In order to show that all Fourier coefficients of $\zeta_{k,\Delta,\rho}(f)$ lie in $i\pi^k\sqrt{\Delta}F_G$ we use Hecke operators. First note that we have
	\begin{align}\label{hecke operator on eta}
	\eta_{k,\Delta,\rho}(f)|T_m = \eta_{k,\Delta,\rho}(m^{2k-1}f|T_m)
	\end{align}
	for any $m \in \N$. This may be checked using the explicit action \eqref{hecke operator} of $T_m$ on the Fourier expansion of $f$ and the action of $T_m$ on $f_{k,D,r,\Delta,\rho}$, which can be computed as in the proof of \cite[Equation~(36)]{zagiereisenstein}. For example, if $p$ is prime with $(p,N) = 1$ we have
	\[
	f_{k,D,r,\Delta,\rho}|T_p = f_{k,p^2D,pr,\Delta,\rho} + p^{k-1}\left(\frac{D\sigma}{p} \right)f_{k,D,r,\Delta,\rho}+p^{2k-1}f_{k,D/p^2,r/D,\Delta,\rho},
	\]
	and there are similar formulas for any Hecke operator $T_m$. We leave the details of the verification of \eqref{hecke operator on eta} to the reader. Now \eqref{hecke operator on eta} implies that
	\begin{align}\label{hecke action zeta}
	\zeta_{k,\Delta,\rho}(f)|T_m - \lambda_m(G)\zeta_{k,\Delta,\rho}(f) = \eta_{k,\Delta,\rho}(f')
	\end{align}
	with the weakly holomorphic modular form $f' = m^{2k-1}f|T_m - \lambda_m(G)f \in M_{\frac{3}{2}-k,\widetilde{\rho}_L}(F_G)$. Again, it follows from \cite{mcgraw} that $f'$ has all coefficients in $F_G$, which by \eqref{fourier expansion canonical modular form Zf} implies that $\eta_{k,\Delta,\rho}(f')$ has all coefficients in $i\pi^k\sqrt{\Delta}F_G$. Using the explicit action of $T_m$ on Fourier expansions and the fact that the first Fourier coefficient of $\zeta_{k,\Delta,\rho}(f)$ vanishes, we obtain by induction from \eqref{hecke action zeta} that all coefficients of $\zeta_{\Delta,\rho}(f)$ lie in $i\pi^k\sqrt{\Delta}F_G$.
\end{proof}

\section{Outlook: Derivatives of $L$-functions of newforms of higher weight}
\label{outlook}

In this last section we explain a possible future application of our results to a non-vanishing criterion for central values of derivatives of $L$-functions of newforms of weight $2k$. We keep the notation from Section~\ref{meromorphic modular forms for heegner divisors} and first recall a non-vanishing criterion for the central $L$-derivatives of newforms of weight $2$ from \cite{bruinieronoheegner}. Let $G \in S_{2}^{\new,\epsilon}(N)$ be a newform of weight $2$ for $\Gamma_0(N)$. The functional equation of the twisted $L$-function $L(G,\chi_\Delta,s)$ implies that $L(G,\chi_\Delta,1) = 0$. Hence, it is natural to consider the (non-)vanishing of $L'(G,\chi_\Delta,s)$ at $s = 1$. The following theorem connects this question to the algebraicity of Fourier coefficients of (the holomorphic part of) harmonic Maass forms.

\begin{theorem}[Theorem~7.6 in \cite{bruinieronoheegner}] Let $G \in S_{2}^{\new,\epsilon}(N)$, let $g \in S_{\frac{3}{2},\overline{\widetilde{\rho}}_L}$ be its Shimura correspondent with coefficients in $F_G$, and let $f \in H_{\frac{1}{2},\rho_L}(F_G)$ be a harmonic Maass form with principal part defined over $F_G$ and $\xi_{\frac{1}{2}}f = \|g\|^{-2}g$. Then the following are equivalent.
\begin{enumerate}\label{bo}
	\item We have $L'(G,\chi_\Delta,1) = 0$.
	\item We have $c_f^+(|\Delta|,\rho) \in F_G$.
\end{enumerate}
\end{theorem}

The proof of this theorem consists of three main steps.
\begin{enumerate}
	\item First, one may show that the canonical differential $\eta_{\Delta,\rho}(f)$ for the (degree $0$) twisted Heegner divisor 
	\[
	y_{\Delta,\rho}(f):=Z_{\Delta,\rho}(f)-\deg(Z_{\Delta,\rho}(f))\cdot\infty
	\]
	is defined over a number field if and only if $c_f^+(|\Delta|,\rho) \in F_G$. This is done by constructing $\eta_{\Delta,\rho}(f)$ as a regularized theta lift and using Hecke operators.
	\item Secondly, by transcendence results for differentials of the third kind due to Scholl \cite{scholl}, Waldschmidt \cite{waldschmidt}, and W\"ustholz \cite{wuest}, we have that the canonical differential $\eta_{\Delta,\rho}(f)$ is defined over a number field if and only if some integral multiple of the Heegner divisor $y_{\Delta,\rho}(f)$ is a principal divisor. This is in turn equivalent to saying that the N\'{e}ron-Tate height of $y_{\Delta,\rho}(f)$ vanishes.
	\item Lastly, by the Gross-Zagier Theorem \cite{gz} the N\'{e}ron-Tate height of $y_{\Delta,\rho}(f)$ is a multiple of $L'(G,\chi_\Delta,1)$, which concludes the proof of the above theorem.
\end{enumerate}

It would be interesting to extend Theorem~\ref{bo} to newforms of higher weight $2k$. We may speculate that the vanishing of $L'(G,\chi_\Delta,k)$ for a newform $G \in S_{2k}^{\new,\epsilon}(N)$ is equivalent to the algebraicity of the coefficient $c_f^+(|\Delta|,\rho)$, where $f \in H_{\frac{3}{2}-k,\widetilde{\rho}_L}(F_G)$ is a harmonic Maass form whose principal part is defined over $F_G$ and which satisfies $\xi_{\frac{3}{2}-k}(f) = \|g\|^{-2}g$, with $g \in S_{\frac{1}{2}+k,\overline{\widetilde{\rho}}_L}^{\new}$ being the Shimura correspondent of $G$. 

Our Theorem~\ref{canonical modular form heegner level N} generalizes step (1) in the above proof sketch to higher weight $2k$. Moreover, step (3) is (essentially) taken care of by Zhang's generalization of the Gross-Zagier formula to newforms of higher weight \cite{zhang}. Here the N\'{e}ron-Tate height on the Jacobian of $X_0(N)$ has to be replaced with a height pairing on the Kuga-Sato variety $W^{2k-2}$ of dimension $2k-1$ over $X_0(N)$, and the Heegner divisor $y_{\Delta,\rho}(f)$ has to be replaced with its corresponding \emph{Heegner cycle} in the Chow group of codimension $k$ cycles on $W^{2k-2}$. Concerning step (2) of the above proof sketch, it remains to show that the Fourier coefficients of the canonical meromorphic modular form $\eta_{k,\Delta,\rho}(f)$ are contained in $i\pi^k\sqrt{\Delta}F_G$ if and only if the Heegner cycle corresponding to $y_{\Delta,\rho}(f)$ vanishes in the Chow group. 

A possible path to prove one direction of this claim is laid out in the thesis of Mellit \cite{mellit}. Given a principal divisor in the Chow group of codimension $k$ cycles in $W^{2k-2}$, Mellit constructs a meromorphic modular form with algebraic Fourier coefficients. We expect that Mellit's meromorphic modular form attached to the Heegner cycle associated with $y_{\Delta,\rho}(f)$ (assuming that it is principal) is a multiple of $\eta_{k,\Delta,\rho}(f)$.

The converse direction seems to be more difficult. Assuming that the coefficients of $\eta_{k,\Delta,\rho}(f)$ are contained in $i\pi^k\sqrt{\Delta}F_G$, we would need to construct a rational function $F$ on a suitable codimension $k-1$ cycle $U$ in the Kuga-Sato variety $W^{2k-2}$ such that the divisor of $F$ is the Heegner cycle corresponding to $y_{\Delta,\rho}(f)$. For $k = 1$, the cycle $U$ is just the whole modular curve $X_0(N)$, and it was shown in \cite{bruinieronoheegner} that the rational function $F$ with divisor $y_{\Delta,\rho}(f)$ can be constructed as the Borcherds product attached to $f$, but for $k > 1$ it is not clear how to construct $U$ and $F$. We plan to come back to this problem in the future.

\bibliographystyle{alpha}
\bibliography{bib.bib}

\end{document}